\documentclass{article}[12pt]
\usepackage{amssymb}
\usepackage[latin1]{inputenc}
\usepackage{amsthm}
\usepackage{latexsym}
\usepackage{amsmath}
\usepackage{graphicx} 
\newcommand{\beqa}{\begin{eqnarray}}
\newcommand{\eeqa}{\end{eqnarray}}
\newcommand{\N}{\mathbb N}
\newcommand{\Z}{\mathbb Z}

\newcommand{\C}{\mathbb C}
\newcommand{\X}{\mathrm X}

\newcommand{\U}{\mathrm U}
\newcommand{\CC}{\mathrm C}

\newcommand{\V}{\mathrm V}

\newcommand{\G}{\mathrm G}
\newcommand{\A}{\mathrm A}

\newcommand{\fimdemo}{\begin{flushright} $\blacksquare$ \end{flushright}}
\newtheorem{teorema}{Theorem}[section]
\newtheorem{corolario}[teorema]{Corolary}

\newtheorem{definicao}[teorema]{Definition}
\newtheorem{prop}[teorema]{Proposition}
\newtheorem{observacao}[teorema]{Remark}
\newtheorem{exemplo}[teorema]{Example}
\newtheorem{exemplos}[teorema]{Examples}

\newcommand{\bethm}{\vspace{0.4\baselineskip}\begin{minipage}{0.9\textwidth}}
\newcommand{\enthm}{\end{minipage}\vspace{0.4\baselineskip}}

\def \and {\hbox {,\quad and \quad }}
\def \for #1{,\quad \forall \,#1}

\begin{document}

\thispagestyle{empty}                          

\begin{center}
{\Large{\bf Envelope Algebras of Partial Actions as Groupoid C*-Algebras}} \\
\vspace{5mm}
{\large Ruy Exel\footnote{* partially supported by CNPq, Brazil}, Thierry Giordano\footnote{** partially supported by a grant from NSERC, Canada}, Daniel Gonçalves}\\
\end{center}  
\vspace{8mm}

\abstract We describe the envelope C*-algebra associated to a partial action of a countable discrete group on a locally compact space as a groupoid C*-algebra (more precisely as a C*-algebra from an equivalence relation) and we use our approach to show that, for a large class of partial actions of $\Z$ on the Cantor set, the envelope C*-algebra is an AF-algebra. We also completely characterize partial actions of a countable discrete group on a compact space such that the envelope action acts in a Hausdorff space.

\section{Introduction}

The concept of partial actions was introduced in \cite{Exel1} and \cite{Mc} and it has been a very important tool in C*-algebras and dynamical systems ever since. As the name suggests, partial actions generalize the notion of an action in a C*-algebra or in a topological space. The problem of deciding whether or not a given partial action is the restriction of some global action (called envelope action) was studied by F. Abadie in \cite{FAbadie}, where, among other things, he shows that the cross product of the envelope C*-algebra by the envelope action is Morita-Rieffel equivalent (previously known as strongly Morita
equivalent) to the partial cross product.

In this paper we are interested in partial actions of a countable discrete group $\G$ (in particular of $\Z$) on a locally compact space $\X$ (particularly the Cantor set, that is, a compact, totally disconnected, with no isolated points, metric space). It is well known that there is a correspondence between partial actions on a Hausdorff locally compact space $\X$ and the partial actions on the C*-algebra $\CC_0(\X)$ (see proposition 3.3 of \cite{Exel1} or \cite{teseMest} for example). 
In \cite{FAbadie}, it is shown that a partial action on a topological space always has an envelope action, which may {\bf not} act on a Hausdorff space (the odometer partial action for example). When the envelope space is Hausdorff the notion of the envelope action in the category of C*-algebras is a rather natural one, but when the envelope space in non Hausdorff the notion of the envelope action has to be reformulated with the use of C*-ternary rings and the introduction of the notion of strong Morita equivalence between partial actions. Although it seems that this approach can not be avoided in general, in the case of a partial action of a countable discrete group on a locally compact space we give a description of the envelope C*-algebra as a C*-algebra from an equivalence relation (viewed as a groupoid). Our approach has the advantage of working for either Hausdorff or non Hausdorff envelope spaces. We also use our description of the envelope algebra to show that it is an AF-algebra, provided we have a partial action of $\Z$  on the Cantor set with some mild assumptions, namely that it arises as a "restriction" of a global action (we should warn the reader that we use the word restriction here with a slight different meaning then what usually appears in the literature). 

The paper is structured as follows. In section 2 we make a quick review of the necessary notions on partial actions and envelope actions. We completely characterize the partial actions of a countable discrete group on a compact space such that the envelope space is Hausdorff in section 3 and finally in section 4 we describe the envelope algebra as a groupoid C*-algebra and show it is an AF-algebra under the assumptions mentioned above.


\section{Partial Actions and Envelope Actions}

\begin{definicao} A partial action of a group $\G$ on a set $\Omega$ is a pair $\theta=(\{\Delta_{t}\}_{t\in \G},\{h_{t}\}_{t\in \G})$, where for each $t\in \G$, $\Delta_{t}$ is a subset of $\Omega$ and $h_{t}:\Delta_{t^{-1}} \rightarrow \Delta_{t}$ is a bijection such that:
\begin{enumerate}
	\item $\Delta_{e} = \Omega$ and $h_{e}$ is the identity in $\Omega$;
	\item $h_{t}(\Delta_{t^{-1}} \cap \Delta_{s})=\Delta_{t} \cap \Delta_{ts}$;
	\item  $h_{t}(h_{s}(x))=h_{ts}(x),$ $x \in \Delta_{s^{-1}} \cap \Delta_{s^{-1} t^{-1}}.$
\end{enumerate}

If $\Omega$ is a topological space, we also require that each $\Delta_{t}$ is an open subset of $\Omega$ and that each $h_{t}$ is a homeomorphism of $\Delta_{t^{-1}}$ onto $\Delta_{t}$. 

Analogously, a pair $\theta = (\{ D_{t} \}_{t \in \G} , \{ \alpha_{t} \}_{t \in \G} )$ is a partial action of $\G$ on a C*-algebra $\A$ if each $D_{t}$ is a closed two sided ideal and each $\alpha_{t}$ is a *-isomorphism of $D_{t^{-1}}$ onto $D_{t}$.
\end{definicao}

Since we are very interested in partial actions of $\Z$, below we give the most important example of such partial actions.

\begin{exemplo}\label{theimportantexample} Let $\X$ be a locally compact space, $\U$ and $\V$ open subsets of $\X$ and $h$ a homeomorphism from $\U$ to $\V$. Let $\X_{-n}= \text{dom} (h^n)$ and $h_n: X_{-n} \rightarrow X_n$ be defined by $h^n$, for $n\in \Z$. Then $\theta=(\{\X_{n}\}_{n\in \Z},\{h_{n}\}_{n\in \Z})$ is a partial action of $\Z$.
\end{exemplo}
\proof

See \cite{Exel1} or \cite{teseMest}.

\begin{exemplo} The Odometer: Let $\X=\{0,1\}^{\infty}=\displaystyle \prod_{\N}\{0,1\}$. Let $\max=1^\infty$ (sequence of all 1´s), $min=0^\infty$ (sequence of all 0´s), $\X_{-1}=\X-\{max\}$, $\X_1=\X-\{min\}$ and $h:\X_{-1} \rightarrow \X_{1}$ be addition of 1 with carryover to the right. Then $\theta=(\{\X_{n}\}_{n\in \Z},\{h_{n}\}_{n\in \Z})$, where $X_{-n}= \text{dom}( h^n)$, is a topological partial action. 
\end{exemplo}

\begin{observacao} With $D_t=\{f\in \CC_0(\X):f|_{\X_t^c}=0\}$, where $\X_t^c$ means the complement of $\X_t$ in $X$, and $\alpha_t:D_{t^{-1}} \rightarrow D_t$ defined by $\alpha_t(f)=f\circ h_t^{-1}$, we have a partial action on the C*-algebra $\CC(\X)$.
\end{observacao}

We now recall the definition of the envelope action in the topological sense.

\begin{definicao} Let $\theta=(\{\X_{t}\}_{t\in \G},\{h_{t}\}_{t\in \G})$ be a partial action. The envelope space, denoted by $\X^e$, is the topological quotient space $\frac{(\G \times \X)}{\sim}$, where $\sim$ is the equivalence relation given by $$(r,x)\sim (s,y) \Leftrightarrow x\in \X_{r^{-1}s} \text{ and } h_{s^{-1}r}(x)=y.$$ The envelope action, denoted by $h^e$, is the action induced in $\X^e$ by the action $h^e_s(t,x)\mapsto(st,x)$.
\end{definicao}

Given a locally compact space $\X$, the definition of the envelope
action and space in the C*-algebraic sense is motivated by the 1-1
relation between partial actions on $\X$ and partial actions on
$\CC_0(\X)$, (see \cite{teseMest} for a proof of this
relation). Basically, if a partial action $\theta=(\{\X_{t}\}_{t\in
\G},\{h_{t}\}_{t\in \G})$ on $\X$ has 
a Hausdorff envelope space $\X^e$ and envelope action $h^e$, then $\CC_0(\X^e)$ is the envelope C*-algebra and $\alpha^e(f)=f\circ (h^e)^{-1}$ is the induced global action associated to the partial action in $\CC_0(\X)$. In \cite{FAbadie}, it is proved that when $\X^e$ is Hausdorff, the partial cross product of $\CC_0(\X)$ by the partial action is 
  Morita-Rieffel
  equivalent to the global cross product of the envelope algebra $\CC_0(\X^e)$ by the envelope action. The problem arises when $\X^e$ is {\bf not} Hausdorff (as for example in the Odometer partial action). The result mentioned above is not valid anymore, as we may have very few continuous functions in $\CC_0(\X^e)$. Abadie, in \cite{FAbadie}, goes around this problem by making use of C*-ternary rings and  introducing the notion of Morita equivalence between partial actions. 
Below we completely characterize partial actions of a countable discrete group on a compact space for which the envelope space is Hausdorff.  Throughout the rest of the paper $\G$ will denote a countable discrete group.

\section{Partial Actions of $\G$ on a compact space such that the envelope space is Hausdorff.}

In \cite{FAbadie} it is shown that a partial action has a Hausdorff envelope space if and only if the graph of the action is closed. Below we give a concrete characterization of partial actions of a countable discrete group on compact spaces for which the envelope space is Hausdorff.

Let $\X$ be a compact set, $\{\X_t,\alpha_t\}$ a partial action of $\G$ on $\X$ and $(\X^e,\alpha^e)$ the envelope space and action respectively. Recall that $\X^e$ is the quotient of $\G \times \X$ by the equivalence relation $(r,x)\sim (s,y) \Leftrightarrow x\in \X_{r^{-1}s} \text{ and } \alpha_{s^{-1}r}(x)=y$, with the quotient topology. We denote the equivalence class of $(n,x)$ in $\X_e$ by $[n,x]$.

\begin{prop}\label{EnvSpcHausdPAclopen} Let $\G$ be a countable discrete group with unit $e$. Then $\X^e$ is Hausdorff if and only if the partial action $\{\X_t,\alpha_t\}$ acts in clopen subsets of $\X$, that is, if and only if $\X_t$ is clopen for each $t\in \G$.
\end{prop}
\proof

First assume that $\X^e$ is Hausdorff. We will show that each $\X_t$ is closed (it is already open by the definition of a partial action).

Suppose there exists $t\in \G$ such that $\X_t$ is not closed (we will show that this implies that $\X^e$ is non-Hausdorff).

Since $\X_t$ is not closed, there exists a sequence $(x_k)_{k\in \N}$ such that $x_k \in X_t$ for all $k$ and such that $x_k \rightarrow x$, where $x \notin X_t$. By the compactness of $\X$, $(\alpha_{t^{-1}}(x_k))_{k\in \N}$ has a converging subsequence and we may pass to this subsequence. We may therefore assume that there exists a sequence $(x_k)_{k\in \N}$ in $X_t$ such that $x_k \rightarrow x$, where $x \notin X_t$ and such that $(\alpha_{t^{-1}}(x_k))_{k\in \N}$ converges to a point $y\in \X$.

We now have that the points $[t^{-1},x]$ and $[0,y]$ can not be separated. Let´s see why:

Suppose that $U$ and $V$ are open, $[t^{-1},x] \in U$ and $[0,y]\in V$. Remember that $U$ is open if and only if $q^{-1}(U)$ is open, where $q$ is the quotient map.

Well, since $x_k \rightarrow x$ and $\alpha_{t^{-1}}(x_k) \rightarrow y$, there exists $N \in \N$ such that $(t^{-1},x_k) \in q^{-1}(U)$ and $(e,\alpha_{t^{-1}}(x_k)) \in q^{-1}(V)$ for all $k> N$. Now notice that $(t^{-1},x_k)\sim (e,\alpha_{t^{-1}}(x_k))$ and hence $[t^{-1},x_k]=[e,\alpha_{t^{-1}}(x_k)]$ and $U\cap V \neq \emptyset$.

\vspace{1pc}

For the converse, we may now assume that each $\X_t$ is a clopen subset of $\X$.

Let $[r,x] \neq [s,y]$ in $\X^e$. So $(r,x)$ is not equivalent to $(s,y)$. We have two possibilities:

- If $x\notin \X_{r^{-1}s}$ then there exists $V_x$ such that $x\in V_x$ and $V_x \cap X_{r^{-1}s} = \emptyset$ (since $\X_{r^{-1}s}$ is clopen).

We may assume $r\neq s$ (if $r=s$ then $x\neq y$ and we find the desired neighborhoods using the fact that $\X$ is Hausdorff).

Take $V=(r,V_x)$ and $U=(s,X)$. Then $i_r(V)$ and $i_s(U)$ are the desired open sets. (Where $i_r(x)=q(r,x)$, $i_s(x)=q(s,x)$ and $q$ is the quotient map)(Also notice that $i_t$ is an open map and the proof is done analogous to what is done in \cite{FAbadie} for the map $i$ in theorem 2.5).

- If $x\in \X_{r^{-1}s}$ then $\alpha_{s^{-1}r}(x)\neq y$.

Let $z=\alpha_{s^{-1}r}(x).$ We have that $z\neq y$. Since $\X$ is Hausdorff, there exists open sets $U_z$ and $U_y$ such that $U_z \cap U_y=\emptyset$, $z\in U_z$ and $y\in U_y$. Take $V_x = \alpha^{-1}_{s^{-1}r}(U_z\cap\X_{s^{-1}r})$, which is an open set. Then $i_r(V_x)$ and $i_s(U_y)$ have the desired properties.

\fimdemo

\begin{observacao} Dokuchaev and Exel have a result in the more general context of partial actions on associative algebras, see theorem 4.5 of \cite{DE}, that implies the proposition above. Still we believe our proof above helps to give the reader a feeling for the space $\X^e$.
\end{observacao}

\begin{prop} Let $\X$ be the Cantor set, $\G$ a countable discrete group and $\{\X_t,\alpha_t\}$ a partial action of $\G$ on $\X$ such that $\X_t$ is clopen for all $t\in \G$. Then the envelope space $\X^e$ is a locally compact Cantor set.
\end{prop}
\proof 
We need to show that $\X^e$ is Hausdorff, locally compact, has a countable basis of clopen sets and has no isolated points (this is equivalent to the characterization of the Cantor set we gave before).

Before we proceed, we note that for each $t\in \G$ the function $i_t(x)=q(t,x)$ is a continuous, open and closed map (we already know that $i_t$ is continuous and open by an argument similar to what is done in \cite{FAbadie} in theorem 2.5). To see that it is a closed map, let $F$ be closed in $\X$. Then $F$ is compact and hence $i_t(F)$ is compact. Since $\X^e$ is Hausdorff we have that $i_t(F)$ is closed.

Notice that $\X^e$ is Hausdorff by proposition \ref{EnvSpcHausdPAclopen}. To prove that $\X^e$ is locally compact, let $[(r,x)]$ in $\X^e$. Since $\X$ is compact, there exists a compact neighborhood, $U_x$, of $x$ in $\X$. But then $i_r(U_x)$ is a compact neighborhood of $[(r,x)]$ in $\X^e$.

Now, let $\{U_n\}_{n\in \G}$ be a countable basis of clopen sets of $\X$. Then $\{i_t(U_n)\}_{n,t\in \G}$ is a countable basis of clopen subsets of $\X^e$.

Finally, we have that $\X^e$ has no isolated points, since if $[(t,x)] \in \X^e$ and $V$ is an open set that contains $[(t,x)]$ then $(t,x)\in q^{-1}(V)$ (we may assume that $q^{-1}(V)$ is of the form $(t,U)$, where $U$ is open in $\X$). So there exists $(y,t) \neq (x,t)$ such that $(y,t) \in q^{-1}(V)$ and hence $[(y,t)] \in V$ and $[(y,t)]\neq[(x,t)]$.

\fimdemo

With the above propositions we completely characterized the envelope actions of partial actions of $\G$ acting on clopen subsets of the Cantor set. The cross product of their envelope C*-algebra by the envelope action is 
  Morita-Rieffel
  equivalent to the partial cross product, see \cite{FAbadie}. The problem is that most of the interesting examples, including the famous odometer (or adding machine), do not satisfy the conditions of the propositions above. Namely they fail to be partial actions on clopen sets. We note here that the majority of examples from partial actions arises as in example \ref{theimportantexample}. In the next section we show how to deal with these examples in a different (and we believe easier) way from what was done in \cite{FAbadie}. As a consequence of our approach we show that the envelope C*-algebra associated to a partial action of $\Z$ on the Cantor set, as in example \ref{theimportantexample}, is an AF-algebra.

\section{The envelope C*-algebra as a groupoid C*-algebra.}

In this section, we start by showing that the envelope C*-algebra associated to a partial action of a countable discrete group $\G$ on a locally compact space can be seen as a C*-algebra of an equivalence relation (seen as a groupoid in the usual way). Before we proceed we need to introduce the notion of core subalgebras, which will be used in our proof that the envelope algebra can be realized as a groupoid C*-algebra.
 

\begin{definicao} Let $A$ be a C*-algebra and let $B\subseteq A$ be a (not necessarily
closed) *-subalgebra.  We shall say that $B$ is a core
subalgebra of $A$ when every representation\footnote{By a representation of a
*-algebra $B$ we mean a multiplicative, *-preserving, linear map
$\pi:B\to{\cal B}(H)$, where $H$ is a Hilbert space.}
  of $B$ is continuous relative to the norm induced from $A$.
\end{definicao}  

Assuming that $B$ is a core subalgebra of $A$, and given a
representation $\pi$ of $B$, we may therefore extend $\pi$ to a
representation $\bar\pi$ of $\bar B$ (the closure of $B$ in $A$).
Since $\bar B$ is a C*-algebra we have by \cite{Dixmier} that
$\bar\pi$ is necessarily contractive.  Therefore we have:

\begin{prop} $B$ is a core subalgebra of $A$ if and only if every
representation of $B$ is contractive.
\end{prop}

\begin{exemplos}\label{Exemplos}
\end{exemplos}
\begin{enumerate}\renewcommand{\labelenumi}{(\roman{enumi})}
	\item Every closed *-subalgebra of a C*-algebra is a core
subalgebra by \cite{Dixmier}.

  \item\label{Exii} Let $B$ be a *-subalgebra of a C*-algebra $A$, such that
$B=\bigcup_{i\in I}B_i$, where each $B_i$ is a core subalgebra of
$A$. Then $B$ is a core subalgebra of $A$.  This is because every
representation of $B$ must be contractive on each $B_i$.

  \item If $\X$ is a locally compact space then $\CC_c(\X)$ is a core subalgebra
of $\CC_0(\X)$. This follows from the fact that $\CC_c(\X)$ is the union of
the closed *-subalgebras $\CC_0(\U)$, where $\U$ ranges in the collection
of all relatively compact open subsets of $\X$.

  \item\label{Exiv} Let ${\cal G}$ be a groupoid satisfying the hypotheses of
corollary 1.22 of \cite{Renault}. Then $\CC_c({\cal G})$ is a core subalgebra of
$\CC^*({\cal G})$ by the same corollary of \cite{Renault}.

  \item Let $B$ be a *-algebra such that
 \begin{equation}\renewcommand{\theequation}{\ref{Exemplos}.\arabic{equation}}
   \left|\left\|b\right\|\right| := \sup \|\pi(b)\|<\infty, \forall \ b\in B,
   \label{eqsupfinito}
 \end{equation}
  where the supremum is taken over the collection of all
representations $\pi$ of $B$.  Then one may define the envelopping
C*-algebra, $\CC^*(B)$, by moding out the elements $b$ such that $\left|\left\|b \right\|\right|=0$, and completing under $\left|\left\|\cdot \right\|\right|$.
  The image of $B$ within $\CC^*(B)$ is therefore a dense core
subalgebra.
  \item If $\G$ is a discrete group then the complex group algebra $\C
\G$ is a core subalgebra of the full group C*-algebra of $\G$.
  \item If $\G$ is non-amenable then $\C \G$ will {\bf not} be a
core subalgebra of the reduced group C*-algebra $\CC^*_r(\G)$, since some
representations of $\C \G$ will extend to a representation of $\CC^*(\G)$
which does not vanish on the kernel of the left regular
representation.

\end{enumerate}

\bigskip
  Let $B$ be a core subalgebra of a C*-algebra $A$. It is then evident
that $B$ satisfies \ref{Exemplos}.\ref{eqsupfinito} and moreover that $\left|\left\|b \right\|\right|
=\|b\|$, where the right hand side refers to the norm of $b$ computed
as an element of $A$.  Supposing that $B$ is dense in $A$, it follows
that $A$ is isomorphic to the envelopping C*-algebra $\CC^*(B)$.  From
this one immediately has:

\begin{prop}\label{isomcore}
  Suppose that $A_1$ and $A_2$ are C*-algebras, and that $B_i$ is a
dense core subalgebra of $A_i$, for $i=1,2$.  If $B_1$ and $B_2$ are
isomorphic as *-algebras, then $A_1$ and $A_2$ are isometrically
*-isomorphic.
\end{prop}

\begin{teorema}\label{teor45}
  Let $A$ be a C*-algebra and let $\{p_i\}_i$ be a family of mutually
orthogonal projections in 
  the multiplier algebra of $A$,  here denoted as $M(A)$. Also let $B$ be
a *-subalgebra of $A$ such that 
  $$
  B\subseteq\bigoplus_{i, j\in I} B\cap (p_iA p_j), 
  $$
  and such that $B\cap (p_iA p_i)$ is a core subalgebra for every
$i\in I$.  Then $B$ is a core subalgebra.
\end{teorema}

\proof  
Given $b\in B$,  by hypothesis we have that 
  $$
  b=\sum_{k, l\in I}a_{kl},
  $$
  where the nonzero summands are finite and each $a_{kl}\in B\cap
(p_kA p_l)$.  We therefore have for all $i,j\in I$, that
  $$
  p_ibp_j =  \sum_{k, l\in I}p_ia_{kl}p_j = a_{ij}, 
  $$
  from where we see that $a_{ij} = p_ibp_j$.  From now on we will
adopt the notation
  $$
  b_{ij}:= p_ibp_j
  \ \ \ \forall \ b\in A,\ \ \ \forall \  i,j\in I,
  $$
  and hence we have for every $b\in B$ that $b_{ij}\in B\cap (p_iA
p_j)$, while $b=\sum_{i, j\in I}b_{ij}$, a sum with finitely many
nonzero terms.

For each finite set of indices $F\subseteq I$, let 
  $$
  B_F = \bigoplus_{i, j\in F} B\cap (p_iA p_j).
  $$
  It is easy to see that $B_F$ is a *-subalgebra of $A$ and we claim
that it is a core subalgebra.  In fact, given a representation $\pi$
of $B_F$, we have for all $i,j\in F$, and all $b_{ij}\in B\cap (p_iA
p_j)$ that
  $$
  \|\pi(b_{ij})\|^2 =
  \|\pi(b_{ij}b_{ij}^*)\| \leq
  \|b_{ij}b_{ij}^*\| =
  \|b_{ij}\|^2,
  $$
  where the crucial second step follows from the fact that
$b_{ij}b_{ij}^*\in B\cap (p_iA p_i)$, the latter being a core
subalgebra by hypothesis.  Given any $b\in B_F$, we then have
  $$
  \|\pi(b)\| \leq
  \sum_{i, j\in F}\|\pi(b_{ij})\| \leq
  \sum_{i, j\in F}\|b_{ij}\| \leq
  |F|^2\|b\|.
  $$
  This proves that $\pi$ is bounded and hence that $B_F$ is a core
subalgebra, as claimed.

Now observe that $B = \bigcup_FB_F$, where $F$ ranges in the
collection of all finite subsets of $F$, so the conclusion follows
from (\ref{Exemplos}.ii).
\fimdemo

We can now focus again in realizing the envelope algebra as a groupoid
C*-algebra. 

  For this we fix,  as before, a partial action of the discrete group
$G$ on a locally compact space $X$.

Recall that the envelope space is the quotient of $\G \times \X$ by the equivalence relation $(r,x)\sim (s,y) \Leftrightarrow x\in \X_{r^{-1}s} \text{ and } \alpha_{s^{-1}r}(x)=y$. So instead of considering this quotient, we will consider the equivalence relation $R\subseteq \G \times \X \times \G \times \X$ given by $$(r,x)\sim (s,y) \Leftrightarrow x\in \X_{r^{-1}s} \text{ and } \alpha_{s^{-1}r}(x)=y,$$ with the product topology.

Notice that a neighborhood base for $z=\left( t,x,s,y \right)$, in $R$, is formed
by neighborhoods of the form $$ \U_{txs}=\{(t,x,s,\alpha_{s^{-1}t}(x)):x\in \U_x\subseteq X_{t^{-1}s} \text{,where $\U_x$ is open} \}.$$

Before we can consider the groupoid C*-algebra of this equivalence relation
  we will show
  that $R$ with this topology is étale, that is, we have to show that $R$ can be equipped with two maps, called range and source, defined by $r\left( t,x,s,y \right)=(t,x)$ and $s\left( t,x,s,y \right)=(s,y)$ and such that $R$ is $\sigma$-compact, $\Delta = \{(t,x,t,x) \in R : (t,x) \in \Z\times \X \}$ is an open subset of $R$ and for all $\left( t,x,s,y \right) \in R$, there exists a neighborhood $\U$ of $\left( t,x,s,y \right)$ in $R$, such that $r$ restricted to $\U$ and $s$ restricted to $\U$ are homeomorphisms from $\U$ onto open subsets of $\X \times \Z$, see \cite{Putnam3} or \cite{Renault}.

\begin{prop} $R$ is étale. 
\end{prop}
\proof

To see that $R$ is sigma compact, we notice that for each fixed $s$ and $t \in \G$, $\X_{t^{-1}s}$ is a countable union of compact sets and hence each $\U_{txs}$, with $\U_x=\X_{t^{-1}s}$ is a countable union of compact sets.

For $(t,x,t,x) \in \Delta$, we take $\U_{txt}$ with $U_x=\X$ to see that $\Delta$ is open.

Finally, given $\left( t,x,s,y \right) \in R$, it is not hard to see that the range and source map are homeomorphisms, once restricted to $\U_{txs}$, with $\U_x=\X_{t^{-1}s}$.
\fimdemo

We are now able to consider the full groupoid C*-algebra of $R$, which we denote by $C^*(R)$ (see \cite{Renault} or \cite{GP} for details on the groupoid C*-algebra of étale equivalence relations). Next we show that $C^*(R)$ is isomorphic to the Morita envelope algebra defined in \cite{FAbadie}. In order to do so, we quickly remind the reader of the definitions in \cite{FAbadie} (adapted to the case at hand). 

Given a Fell Bundle $B=(B_t)_{t\in \G}$ of a partial action $\{\X_t,h_t\}_{t\in \G}$, 
        (in our case $B_t=\CC_0(\X_t)\delta_t$), 
we consider the linear space, $k_c(B)$, of all continuous functions
$k:\G \times \G \rightarrow B$,  with
  finite
  support and such that 
  $k(r,s) \in B_{rs^{-1}}$. 
  We now equip $k_c(B)$ with the involution $k^*(r,s)=k(s,r)^* \text{, }\forall k\in k_c(B)$, the multiplication $k_1 * k_2 (r,s) = \displaystyle \sum_{t\in \G} k_1(r,t) k_2(t,s) \text{, } \forall k_1, k_2 \in k_c(B)$ and the norm $\|k\| = \left(\sum_{r,s\in\G} \|k(r,s)\|^2 \right)^{\frac{1}{2}}$. The universal C*-algebra of the completion of $k_c(B)$ with respect to the norm above is the envelope algebra, $k(B)$. Finally we notice that there exists a natural action of $\G$ on $k_c(B)$, which can be extended to $k(B)$, given by $\beta_t(k)(r,s)=k(rt,st)$. The pair $(k(B),\beta)$ is the envelope action as in \cite{FAbadie}. We can now prove our main result.

\begin{teorema} Given a partial action $h$ of a countable discrete
group $\G$ on a locally compact space $\X$, the groupoid C*-algebra
$C^*(R)$, as defined above, is isomorphic to the envelope C*-algebra
$k(B)$.
\end{teorema}

\proof
  Initially let us observe that, given any element $(r,x,s,y)\in R$,
one has that $y=h_{s^{-1} r}(x)$.  Therefore the fourth variable ``$y$"
is a function of the first three, and hence we may discard it.  In
more precise terms we have that
  $$
  (r,x,s,y)\mapsto (x,r,s)
  $$
  establishes a one-to-one correspondence from $R$ to the set
  $$
  R' = \{(x,r,s)\in \X\times\Z\times\Z: x\in \X_{r^{-1} s}\}.
  $$

Moreover this correspondence is seen to be a homeomorphism if $R'$ is
viewed as a subspace of the topological product space
$\X\times\Z\times\Z$.

Borrowing the groupoid structure from $R$ we have that $R'$ itself
becomes an \'étale groupoid under the multiplication operation
  $$
  (x,r,s) \cdot (y,t, u) = (x,r,u),
  $$
  defined iff $y =  h_{s^{-1} r}(x)$, and $s=t$, while the
inversion operation is given by
  $$
  (x,r,s)^{-1} = \big(h_{s^{-1} r}(x),s,r\big).
  $$
  Since $R$ and $R'$ are isomorphic topological groupoids, it is
enough to show that $C^*(R')$ and $k(B)$ are isomorphic C*-algebras.
We will derive this result from proposition \ref{isomcore}, by
showing the existence of two isomorphic dense core subalgebras of
$C^*(R')$ and $k(B)$, respectively.

On the one hand recall that $C_c(R')$ is a dense core subalgebra of
$C^*(R')$, as observed in \ref{Exemplos}.iv.  To define the
relevant dense core subalgebra of $k(B)$, recall that $B$ is the Fell
bundle with fibers $B_t = C_0(X_t)\delta_t$.

Put $D_t = C_c(X_t)\delta_t$, so that each $D_t$ is a dense linear
subspace of $B_t$.  Moreover it is easy to see that
  $$
  D_rD_s\subseteq D_{rs} \and D_r^*=D_{r^{-1}}, 
  \eqno{(\dagger)}
  $$
  for all $r,s\in G$.

Denote by $k_c(D)$ the subset of $k_c(B)$ formed by all $k\in k_c(B)$
such that $k(r,s)\in D_{rs^{-1}}$, for all $r$ and $s$.  As a
consequence of $(\dagger)$ one easily proves that $k_c(D)$ is a
*-subalgebra of $k(B)$, which is also easily seen to be dense.

We will next show that $k_c(D)$ is a core subalgebra of $k(B)$ by
using theorem \ref{teor45}.  With this in mind we must first define a
family of projections $\{p_t\}_{t\in G}$ in the multiplier algebra
$M(k(B))$.
  Given $t\in G$, consider the maps
  $$
  L_t, R_t: k_c(B) \to k_c(B),
  $$
  given,  for every $k\in k_c(B)$,  by
  $$
  L_t(k)(r,s) = \begin{cases}
    k(r,s), & \hbox{ if } r=t, \\ 0, & \hbox{ otherwise}
    \end{cases},
  \qquad \forall r,s\in G.
  $$
  and
  $$
  R_t(k)(r,s) = \begin{cases}
    k(r,s), & \hbox{ if } s=t, \\
    0, & \hbox{ otherwise}
    \end{cases},
  \qquad \forall r,s\in G.
  $$
  One may then prove that both $L_t$ and $R_t$ extend continuously to
$k(B)$ giving a multiplier
  $$
  p_t = (L_t,R_t)\in M\big(k(B)\big),
  $$
  which turns out to be self-adjoint and idempotent, thus producing a
family $\{p_t\}_{t\in G}$ of mutually orthogonal projections.

For every $r,s\in G$ one has that $p_rk_c(D)p_s$ consists of all the
$k\in k_c(D)$ which are supported on the singleton $\{(r, s)\}$.  In
particular $p_tk_c(D)p_t = C_c(\X_t)$, which sits within $p_tk(B)p_t =
\CC_0(\X_t)$ as a core subalgebra, as in example \ref{Exemplos}.iii.  One
may then easily prove that theorem \ref{teor45} applies and hence we
deduce that $k_c(D)$ is a core subalgebra of $k(B)$.

We will next prove that $C_c(R')$ and $k_c(D)$ are isomorphic as
*-algebras and hence the result will follow from proposition \ref{isomcore}.
  Given $r, s\in G$, let 
  $$
  R'_{r, s} = R'\cap \big(X\times\{r\}\times\{s\}\big) 
  $$
  or, equivalently,
  $$
  R'_{r, s} =
  \{(x,r,s): x\in \X_{r^{-1} s}\},
  $$
  so $R'_{r,s}$ naturally identifies with $\X_{r^{-1} s}$.  Given $f\in
C_c(R')$, denote by $f_{r,s}$ the restriction of $f$ to $R'_{r, s}$,
seen as an element of $C_c(\X_{r^{-1} s})$.  Since $f$ is compactly
supported, only finitely many $f_{r,s}$ will be nonzero.  Define
$\psi:C_c(R') \to k_c(D)$ by
  $$
  \psi(f)(r,s) = f_{r^{-1},s^{-1}}\delta_{r s^{-1}}
  \for f\in C_c(R')
  \for r,s\in G.
  $$

Observing that $R'$ is the disjoint union of the $R'_{r, s}$, it
should be obvious that $\psi$ is a well defined vector space
isomorphism.  The proof will then be concluded once we show that
$\psi$ is a *-homomorphism.

In order to prove that $\psi(f*g) = \psi(f)\psi(g)$, we may suppose
without loss of generality that $f$ is supported in $R'_{r, s}$ and
that $g$ is supported in $R'_{t, u}$.  When $s\neq t$, the product in
$R'$ of $(x,r,s)$ and $(y,t,u)$ is never defined, so $f*g=0$.

Otherwise, if $s=t$, we have that $f*g$ is supported in $R'_{r,u}$.
Moreover, given $(x,r,u)\in R'_{r,u}$, the only way of writting
$(x,r,u)$ as a product of an element of $R'_{r, s}$ and an element of
$R'_{s, u}$ is
  $$
  (x,r,u) = (x,r,s) \big(h_{s^{-1} r}(x),s,u\big),
  $$
  as long as $x\in X_{r^{-1} s}$.
  Thus 
  $$
  (f*g)(x,r,u) =
  f(x,r,s) g\big(h_{s^{-1} r}(x),s,u\big) =
  f_{r,s}(x) g_{s, u}\big(h_{s^{-1} r}(x)\big).
  $$

On the other hand, since $\psi(f)$ is supported on the singleton
$\{(r^{-1}, s^{-1})\}$, and $\psi(g)$ is supported on $\{(s^{-1},
u^{-1})\}$, we have that $\psi(f)\psi(g)$ is supported on $\{(r^{-1},
u^{-1})\}$, and
$$\begin{array}{rl}
  \big(\psi(f)\psi(g)\big) (r^{-1}, u^{-1})= & 
  \psi(f)(r^{-1}, s^{-1})  \  \psi(g)(s^{-1}, u^{-1}) \\ = &
  (f_{r,s}\delta_{r^{-1} s}) \  (g_{s,u}\delta_{s^{-1} u}) =
  f_{r,s} (g_{s,u}\circ h_{s^{-1} r}) \delta_{r^{-1} u},
  \end{array}$$
  from where it is easily seen that $\psi(f*g) = \psi(f)\psi(g)$.
  We leave it for the reader to prove that $\psi$ preserves the
adjoint operation.  
\fimdemo

\begin{corolario} Let $\alpha$ be the action on $\CC^*(R)$ given by $\alpha_t(f)(r,x,s,y)=f(rt,x,st,y)$. Then $\CC^*(R) \rtimes_{\alpha} \G$ is isomorphic to $k(B) \rtimes_{\beta} \G$, which is strong morita equivalent to the partial cross product $\CC(\X)\rtimes \G$.
\end{corolario}
\proof
It is clear that the actions $\alpha$ and $\beta$ commute and hence the isomorphism follows. The second part is done in \cite{FAbadie}
\fimdemo

We finish the paper showing that for partial actions of $\Z$ on the Cantor set, $\X$, as in \ref{theimportantexample}, with $\X_{-1}\neq \X$, $R$ is an approximately proper equivalence relation and $C_r^*(R)$ (and hence the envelope C*-algebra) is an AF-algebra.

Recall that an equivalence relation is said to be proper when the quotient space is Hausdorff. In \cite{Re2}, Renault defines approximately proper and a approximately finite equivalence relations as below.

\begin{definicao} An equivalence relation $R$, on a locally compact, second countable, Hausdorff space $\X$, is said to be approximately proper if there exists an increasing sequence of proper equivalence relations $\{R_n
\}_{n\in \N}$ such that $R=\bigcup_{n\in \N} R_n$. An approximately proper equivalence relation on a totally disconected space is called an AF equivalece relation. \end{definicao}

\begin{observacao} In \cite{GPS2}, Giordano, Putnam and Skau define an AF equivalence relation as an equivalence relation that can be written as an increasing union of compact open étale subequivalence relations. They also mention that their definition is equivalent to the definition above.
\end{observacao}

To prove that $R$ is approximately proper we will come up with a sequence of partial actions by clopen sets (so that their envelope space is Hausdorff) such that the union of the induced equivalence relations is $R$.

Recall that $R$ is associated with a partial action $\theta=\{\X_{-n},h^n\}$ on $\X$, as in example \ref{theimportantexample}. That is, $h$ is a homeomorphism from $\U$ to $\V$ (where $\U$ is a proper open subset of $\X$), $\X_{-n}= \text{dom} (h^n)$ and $h_n =h^n$.

To create the partial actions, let $\{\U_k\}_{k=0,1,\ldots}$ be an increasing sequence of clopen sets such that their union is $\X_{-1}=\U\neq \X$. For each $\U_k$, denote the partial action by clopen sets obtained by restricting $h$ to $\U_k$ and proceeding as in example \ref{theimportantexample} by $\theta_k=\{\X^k_{-n},h_n\}_{n\in \Z}$, where $\X^k_{-1}=\U_k$ and $h_1$ is $h$ restricted to $\U_k$.    

Now, we consider the sub equivalence relation $R_k\subseteq \Z\times \X \times  \Z\times \X $ given by $(r,x)\sim_k (s,y) \Leftrightarrow x\in \X^k_{r^{-1}s} \text{ and } h_{s^{-1}r}(x)=y$. Since each $R_k$ is associated to a partial action on clopen sets, we have by proposition \ref{EnvSpcHausdPAclopen} that the quotient $\frac{\Z \times \X}{\sim_k}$ is Hausdorff for every $k$. With this set up we can prove that $R$ is approximately proper.

\begin{prop}\label{Gaproper} $R$ is approximately proper.
\end{prop}
\proof

It remains to show that $R=\bigcup_{k\in \N} R_k$ (it is clear that $R_k \subseteq G_{k+1}$ for $k=0,1,\ldots$). It follows promptly that $R_k \subseteq R$ for all $k$. Next we show that $R\subseteq \bigcup_{k\in \N} R_k$.

Let $(r,x,s,y) \in R$ (which happens iff $x\in \X_{r^{-1}s} \text{ and } h_{s^{-1}r}(x)=y$). All we need to do is find a $K$ such that $x\in \X_{r^{-1}s}^K$, since this would imply that $(r,x,s,y) \in R_K$. Now recall that $\X_{-n}= \text{dom}(h^n)$ and assuming that $r^{-1}s\geq 0$ (the case $r^{-1}s \leq 0$ is analogous) we have that $$\X_{r^{-1}s}=\text{dom}(h^{r^{-1}s})= \U \cap h^{-1}(\U)\cap \ldots \cap h^{s^{-1}r+1}(\U).$$ So $x, h(x), h^2(x), \ldots, h^{r^{-1}s+1}$ belong to $\U$ and hence we can find a $K$ such that $x, h(x), h^2(x), \ldots, h^{r^{-1}s+1}$ all belong to the same $\U_K$, since $\U = \bigcup_{k\in \N} \U_k$ with $\U_k \subseteq U_{k+1}$. We conclude that $x\in \U_K \cap h^{-1}(\U_K)\cap \ldots \cap h^{s^{-1}r+1}(\U_K) = \X_{r^{-1}s}^K$ as desired.
\fimdemo

Since we have shown that $R$ is approximately proper, it is natural to consider $R= \cup R_k$ with the inductive limit topology. This approach will allow us to write $C_r^*(R)$ as an inductive limit C*-algebra. But first we need to show that the inductive limit and product topology agree on $R$.

\begin{prop} Let $R = \bigcup_{k\in\Z}R_k $ above. Then the inductive limit topology and the product topology on $R$ are the same.
\end{prop}
\proof

Suppose $\U\neq \emptyset$ is open in the inductive limit topology. Then $\U \cap R_k$ is open for all $k$. Let $(t,x,s,y)\in \U$ and find $K$ such that $(t,x,s,y) \in R_K$. Then $\U \cap R_K$ contain an open neighborhood of $(t,x,s,y)$ of the form $$\{(t,z,s,h_{s^{-1}t}(z)): z \in \U^K \subseteq \X^K_{t^{-1}s} \subseteq \X_{t^{-1}s}\},$$ where $\U^K$ is open in $\X^K_{t^{-1}s}$ and hence open in $\X_{t^{-1}s}$. So $\U$ is open in the product topology.

Now, notice that $$ \U_{txs}\cap R_k=\{(t,x,s,h_{s^{-1}t}(x):x\in \U_x\subseteq X_{t^{-1}s} \text{,where $\U_x$ is open} \}\cap R_k$$ is open in $R_k$ for all $k$ and hence $ \U_{txs}$ is open in the inductive limit topology.
\fimdemo



\begin{corolario} $C_r^*(R)= \displaystyle \lim_{\longrightarrow} C_r^*(R_k)$.
\end{corolario}
\proof
The proof is analogous to what is done in \cite{GP} for C*-algebras from substitution tilings.
\fimdemo

\begin{prop} $C_r^*(R)$ is an AF-algebra.
\end{prop}
\proof

We already know, by proposition \ref{Gaproper}, that $R$ is approximately proper. Also, $R$ is an equivalence relation in $\Z \times \X$ and since $\X$ is the Cantor set it is clear that $R$ is and AF equivalence relation. Then by theorem 3.9 of \cite{GPS2} we have that $R$ is isomorphic to tail equivalence in some Bratteli diagram and by \cite{ExRe} we have that the associated C*-algebra is an AF algebra.

Another way to prove this proposition would be to show that each $R_k$ is an AF equivalence relation, (as defined in \cite{GPS2}), so that $C_r^*(R_k)$ is an AF C*-algebra. Since inductive limits of AF C*-algebras are again AF (via the local characterization of AF algebras) this will yield that $C_r^*(R)$ is also AF. To see that each $R_k$ is AF, notice that $\cup_n R_k^n = R_k$, where $ R_k^n = \{ (r,x,s,y)\in \Z \times \X \times \Z \times \X : |r|,|s|\leq n \} \cap R_k.$    


\fimdemo

\vspace{1.5pc}

R. Exel, Departamento de Matemática, Universidade Federal de Santa Catarina, Florianópolis, 88040-900, Brasil

Email: exel@mtm.ufsc.br

\vspace{0.5pc}
T. Giordano, Department of Mathematics and Statistics, University of Ottawa, Ottawa, K1N 6N5, Canada

Email: giordano@uottawa.ca
\vspace{0.5pc}

D. Goncalves, Departamento de Matemática, Universidade Federal de Santa Catarina, Florianópolis, 88040-900, Brasil

Email: daemig@gmail.com
\end{document}